\documentclass[11pt]{article}
\usepackage[T1]{fontenc}
\usepackage{lmodern,microtype,mathtools,amssymb,amsthm}
\usepackage[margin=1.02in]{geometry}
\usepackage{xcolor}
\usepackage{enumitem,needspace}
\setlist[enumerate,1]{label=\textup{(\roman*)},ref=\roman*,leftmargin=*,topsep=4pt,itemsep=3pt,parsep=0pt}
\numberwithin{equation}{section}
\usepackage[backend=biber,style=numeric,sorting=nyt,doi=false,isbn=false]{biblatex}
\DeclareFieldFormat[article]{title}{\iffieldundef{doi}{\mkbibemph{#1}}{\href{https://doi.org/\thefield{doi}}{\mkbibemph{#1}}}}
\DeclareFieldFormat[unpublished]{title}{\mkbibemph{#1}}
\renewbibmacro*{in:}{}
\usepackage[colorlinks=true,linkcolor=blue!45!black,urlcolor=blue!45!black,citecolor=blue!45!black]{hyperref}
\newtheorem{theorem}{Theorem}[section]
\newtheorem{corollary}[theorem]{Corollary}
\newtheorem{lemma}[theorem]{Lemma}
\newtheorem{proposition}[theorem]{Proposition}
\newtheorem{fact}[theorem]{Fact}
\newtheorem{claim}{Claim}[theorem]
\theoremstyle{definition}
\newtheorem{definition}[theorem]{Definition}
\newtheorem{remark}[theorem]{Remark}
\theoremstyle{plain}
\DeclareMathOperator{\Av}{Av}
\DeclareMathOperator{\Stab}{Stab}
\DeclareMathOperator{\supp}{supp}
\DeclareMathOperator{\tp}{tp}
\newcommand{\res}{\mathord{\upharpoonright}}
\newcommand{\1}{\mathbf{1}}
\newcommand{\U}{\mathfrak C}
\newcommand{\measureSpace}[2]{\mathfrak{M}_{#1}(#2)}
\title{Idempotent generically stable measures in NIP groups}
\author{Artem Chernikov}
\date{}
\begin{document}
\maketitle
\begin{abstract}
Let $\mu$ be a global generically stable Keisler measure on a type-definable group in an NIP theory. Then $\mu$ is idempotent if and only if it is the unique translation-invariant measure on its type-definable stabilizer subgroup. This resolves the NIP case of the Idempotent Measure problem studied by Chernikov, Gannon, and Krupi\'nski~\cite{CG,CGII,CGK}.
\end{abstract}

\section{Introduction}

A regular Borel probability measure on a locally compact group is idempotent under convolution if and only if it is the normalized Haar measure of a compact subgroup, namely its support~\cite{Pym}. The compact case goes back to Kawada and It\^o~\cite{KI} and Wendel~\cite{Wendel}. Rudin~\cite{Rudin} and Cohen~\cite{Cohen} studied idempotent measures on locally compact abelian groups, and Glicksberg~\cite{Glicksberg} treated compact abelian semigroups.

An analogous question for Keisler measures was considered in~\cite{CG}. On a definable group $G$, these are finitely additive probability measures on the definable subsets of $G$, with convolution given by the pushforward of the Morley product under multiplication. In stable groups the global idempotent types are precisely the generic types of connected type-definable subgroups~\cite{Newelski}. It was proved in~\cite[Theorem~5.8]{CG} that in stable theories every idempotent measure is the unique translation-invariant measure on its type-definable stabilizer. In~\cite[Theorem~4.11]{CGII}, it was shown that in NIP theories a global idempotent measure invariant over a small model and under right translation by $G^{00}$ concentrates on its type-definable right stabilizer.
The Idempotent Measure problem~\cite[Problem~3.41]{CGK} asks whether every global idempotent frequency interpretation measure (fim measure) on a type-definable group is generically transitive, or equivalently, whether it concentrates on its right stabilizer~\cite[Proposition~3.37]{CGK}. In NIP theories the fim condition is equivalent to generic stability~\cite[Theorem~3.2]{HPS}. The problem was settled for abelian type-definable groups in arbitrary theories~\cite[Theorem~3.45]{CGK}. For generically stable types the problem was settled in rosy theories and for inp-minimal groups, including dp-minimal groups~\cite[Propositions~2.37 and~2.25]{CGK}.
Later, Gannon, Hoffmann, and Krupi\'nski formulated an automorphism group version of the problem~\cite[Conjecture~(A)]{GHK}, studied further in~\cite[Theorem~3.14]{HR}. A counterexample to this version in the theory of dense linear orders is given in~\cite[Theorem~A.1]{Hoffmann}, and a generalization of the result for measures from abelian to virtually solvable definable groups is given in~\cite[Corollary~4.8]{Hoffmann}.

Here we answer the Idempotent Measure problem for all groups in NIP theories:

\begin{theorem}\label{thm:main}
Let $G$ be a type-definable group in an NIP theory and let $\mu$ be a global generically stable Keisler measure on $G$. Then $\mu$ is idempotent, i.e., $\mu=\mu*\mu$, if and only if it concentrates on its type-definable right stabilizer $H\coloneqq\Stab_{\mathrm r}(\mu)$. In this case $H=\Stab_{\ell}(\mu)$ and $\mu^{-1}=\mu$, where $\mu^{-1}$ denotes the pushforward under inversion. Among global measures concentrating on $H$, $\mu$ is the unique left-$H$-invariant measure and the unique right-$H$-invariant measure, and $H$ has finitely satisfiable generics (fsg).
\end{theorem}

\begin{corollary}\label{cor:types}
Let $G$ be a type-definable group in an NIP theory and let $p$ be a global generically stable idempotent type on $G$. Then $p$ is generically transitive in the sense of~\cite[Definition~2.13]{CGK}. Explicitly, if $G$, its operations, and $p$ are defined over $M$ and $(a,b)\models p^{(2)}\res M$, then $(ab,b)\models p^{(2)}\res M$.
\end{corollary}
\begin{proof}
Apply Theorem~\ref{thm:main} to the $\{0,1\}$-valued measure associated with $p$ and use~\cite[Remark~2.12]{CGK}.
\end{proof}

As in the abelian case~\cite[Theorem~3.45]{CGK}, we use bounded local weight for generically stable measures (Fact~\ref{fact:average}). In the abelian proof one can move any factor of a product to the end. Here the main technical Lemma~\ref{lem:products} shows that the pushforward of $\mu^{(2k+1)}$ under the map taking the overlapping products $x_{i-1}y_i x_i$, $1\leq i\leq k$, and retaining $y_1,\ldots,y_k$ is $\mu^{(2k)}$, even though successive products share $x_i$.
Since fim measures are closed under Morley products in arbitrary theories~\cite[Theorem~1.3]{CGH}, the proof given here for fim idempotent measures uses NIP only in Proposition~\ref{prop:sym}, in the form of bounded alternation for indiscernible measures (Fact~\ref{fact:indiscernible-limit}, i.e., the measure formulation of Ben Yaacov's theorem that Keisler randomization preserves NIP~\cite{BY}).

\section{Generically stable measures, averaging and pushforwards}

We work in a monster model $\U$ of a complete NIP theory $T$. Fix a small model $M\preccurlyeq\U$ and a group $G$ in variables $x$, type-definable over $M$ with operations defined over $M$.

\begin{definition}
For $A\subseteq\U$, let $\measureSpace{x}{A}$ denote the space of Keisler measures in the variables $x$ over $A$. Thus $\measureSpace{x}{\U}$ is the space of global Keisler measures in $x$. Let $\lambda\in\measureSpace{x}{\U}$ and let $f\colon\U^x\to\U^y$ be definable. The \emph{pushforward} $f_*\lambda\in\measureSpace{y}{\U}$ is given by
\[
    (f_*\lambda)(\varphi(y))\coloneqq \lambda(\varphi(f(x)))
\]
for every $\varphi(y)\in L(\U)$~\cite[Definition~3.25]{CGK}. For a subtuple $u$ of $x$, write $\lambda|_u\coloneqq (\pi_u)_*\lambda$, where $\pi_u$ is the coordinate projection onto $u$. For $A\subseteq\U$, write $\lambda\res A\in\measureSpace{x}{A}$ for the restriction of $\lambda$ to formulas with parameters in $A$.
\end{definition}

A measure \emph{concentrates on} a partial type if it gives measure one to each formula in the partial type. A \emph{measure on $G$} is a global measure concentrating on $G$. We denote the set of these measures by $\measureSpace{G}{\U}$.

\begin{definition}
Let $\mu\in\measureSpace{x}{\U}$ be definable over $M$ and let $\nu\in\measureSpace{y}{\U}$. Their \emph{Morley product} $\mu\otimes\nu\in\measureSpace{xy}{\U}$ is given by
\begin{equation}\label{eq:integral}
  (\mu\otimes\nu)(\theta(x,y;d))
    \coloneqq \int_{S_y(N)}\mu(\theta(x,c;d))\,d(\nu\res N)(q),
       \qquad c\models q ,
\end{equation}
where $N$ is a small model containing $M$ and $d$~\cite[Definition~3.3]{CGK}. By definability of $\mu$ over $M$, the value $\mu(\theta(x,c;d))$ depends only on $\tp(c/N)$ and the resulting function $S_y(N)\to[0,1]$ is continuous. We identify $\nu\res N$ with the corresponding regular Borel probability measure on $S_y(N)$. Its support $\supp(\nu\res N)$ consists of the types all of whose formulas have positive measure, and has measure one. Define the Morley powers by $\mu^{(1)}_{x_1}\coloneqq \mu_{x_1}$ and $\mu^{(r+1)}_{x_1,\ldots,x_{r+1}}\coloneqq \mu_{x_{r+1}}\otimes\mu^{(r)}_{x_1,\ldots,x_r}$ for $r\geq1$. Write $\mu^{(\omega)}$ for the unique measure in $(x_i)_{i\geq1}$ satisfying $\mu^{(\omega)}|_{x_1,\ldots,x_r}=\mu^{(r)}_{x_1,\ldots,x_r}$ for every $r\geq1$.

If $\mu,\nu\in\measureSpace{G}{\U}$, their \emph{convolution} $\mu*\nu\in\measureSpace{G}{\U}$ is the pushforward of $\mu\otimes\nu$ under multiplication on $G$ (see Definition~\ref{def:word-map} and Remark~\ref{rem:agreement}). For $g\in G(\U)$ let $\delta_g\in\measureSpace{G}{\U}$ be the point mass at $g$. The right and left stabilizers of $\mu$ are
\[
  \Stab_{\mathrm r}(\mu)\coloneqq \{g\in G(\U):\mu*\delta_g=\mu\},
  \qquad
  \Stab_{\ell}(\mu)\coloneqq \{g\in G(\U):\delta_g*\mu=\mu\}.
\]
Both are type-definable subgroups over $M$~\cite[Fact~3.35]{CGK}.
\end{definition}

\begin{fact}[{\cite[Theorem~3.2]{HPS}}]\label{fact:fim}
Let $\mu\in\measureSpace{x}{\U}$ be $M$-invariant. In an NIP theory, $\mu$ is generically stable over $M$ if and only if it is both definable over $M$ and finitely satisfiable in $M$. This is also equivalent to $\mu$ being fim over $M$, i.e., for every $L$-formula $\varphi(x;y)$ there are $M$-definable sets $\Theta_r\subseteq(\U^x)^r$ with $\mu^{(r)}(\Theta_r)\to1$ such that, for each $\varepsilon>0$ and all sufficiently large $r$,
\[
 \left|\Av(\bar a)(\varphi(x;c))-\mu(\varphi(x;c))\right|<\varepsilon
 \qquad\text{for all }\bar a\in\Theta_r\text{ and }c\in\U^y.
\]
Here $\Av(\bar a)(\psi(x))\coloneqq\frac1r\sum_{i=1}^r\1_{\psi(a_i)}$ for $\bar a=(a_1,\ldots,a_r)$ and $\psi(x)\in L(\U)$.
\end{fact}

\begin{fact}[{\cite[Fact~3.7 and Proposition~3.26]{CGK}}]\label{fact:products}
Let $\mu\in\measureSpace{x}{\U}$ and $\nu\in\measureSpace{y}{\U}$.
\begin{enumerate}
\item\label{part:associativity} If both measures are definable over $M$, then $\mu\otimes\nu$ is definable over $M$~\cite[Lemma~1.5]{HPS} and $(\mu\otimes\nu)\otimes\lambda=\mu\otimes(\nu\otimes\lambda)$ for every $\lambda\in\measureSpace{z}{\U}$.
\item\label{part:pushforwards} Suppose $\mu$ is definable over $M$ and $\nu$ is arbitrary. For maps $f\colon\U^x\to\U^u$ and $g\colon\U^y\to\U^v$ definable with parameters from $\U$, we have $(f\times g)_*(\mu\otimes\nu)=f_*(\mu)\otimes g_*(\nu)$. Indeed, choose $N\supseteq M$ in~\eqref{eq:integral} containing the parameters of $f$, $g$ and the formula being evaluated. Then $f_*\mu$ is definable over $N$, and the identity follows by change of variables.
\item\label{part:fimclosure} If $\mu$ and $\nu$ are fim over $M$, then $\mu\otimes\nu$ is fim over $M$~\cite[Theorem~1.3]{CGH}. An $M$-definable pushforward of a measure fim over $M$ is again fim over $M$.
\item\label{part:symmetry} If $\mu$ is fim over $M$, then $\mu\otimes\nu=\nu\otimes\mu$ for every definable $\nu$. The finite Morley powers of $\mu$ are invariant under permutations of the coordinates (since $\mu$ commutes with itself).
\end{enumerate}
\end{fact}

\begin{remark}\label{rem:blocks}
Let $\mu\in\measureSpace{x}{\U}$ be generically stable over $M$. 
\begin{enumerate}
\item\label{part:blocks-general} Fix positive integers $m,s,n_1,\ldots,n_s$. For $1\le j\le s$, let $\iota_j\colon\{1,\ldots,n_j\}\hookrightarrow\{1,\ldots,m\}$ be injections with pairwise disjoint images. Define $\pi_j(a_1,\ldots,a_m)\coloneqq(a_{\iota_j(1)},\ldots,a_{\iota_j(n_j)})$. Given maps $f_j\colon(\U^x)^{n_j}\to\U^{y_j}$ definable with parameters from $\U$, put $f\coloneqq(f_1\circ\pi_1,\ldots,f_s\circ\pi_s)$. Then
\[
       f_*(\mu^{(m)})=(f_1)_*(\mu^{(n_1)})\otimes\cdots
                              \otimes(f_s)_*(\mu^{(n_s)}).
\]
(Indeed, by the symmetry in Fact~\ref{fact:products}\textup{(\ref{part:symmetry})}, the projection $(\pi_1,\ldots,\pi_s)$ sends $\mu^{(m)}$ to $\mu^{(n_1+\cdots+n_s)}$. Parts~\textup{(\ref{part:associativity})} and~\textup{(\ref{part:pushforwards})} of that fact let us group the factors into blocks of lengths $n_1,\ldots,n_s$ and apply the pushforward identity successively to $f_1,\ldots,f_s$.)
\item\label{part:blocks-one} Suppose $\mu\in\measureSpace{G}{\U}$ and $\mu*\mu=\mu$. For every $n\geq1$, the ordered product $m_n(x_1,\ldots,x_n)\coloneqq x_1\cdots x_n$ satisfies $(m_n)_*(\mu^{(n)})=\mu$. (Indeed, by part~\textup{(\ref{part:blocks-general})} and the definition of convolution, $(m_{n+1})_*\mu^{(n+1)}=\big((m_n)_*\mu^{(n)}\big)*\mu$. The claim follows by induction from idempotence.)
\item\label{part:blocks-several} Under the assumptions of part~\textup{(\ref{part:blocks-one})}, take $f_j\coloneqq m_{n_j}$ in part~\textup{(\ref{part:blocks-general})} for every $j$. Then $f_*(\mu^{(m)})=\mu^{(s)}$, by the two preceding parts.
\end{enumerate}
\end{remark}

\begin{lemma}\label{lem:fibres}
Let $\mu\in\measureSpace{x}{\U}$ be definable over $M$, let $\nu\in\measureSpace{y}{\U}$, and let $\lambda\in\measureSpace{z}{\U}$ be fim over $M$. Let $f\colon\U^x\times\U^y\to\U^z$ be $M$-definable. Put $f_b(x)\coloneqq f(x,b)$ and $g(x,y)\coloneqq(f(x,y),y)$. Let $P$ consist of the realizations of the types in $\supp(\nu\res M)$. Then
\begin{equation}\label{eq:fibres}
 g_*(\mu\otimes\nu)=\lambda\otimes\nu
 \quad\Longleftrightarrow\quad
 (f_b)_*\mu=\lambda\ \text{for all }b\in P.
\end{equation}
\end{lemma}
\begin{proof}
For the forward implication see~\cite[Lemma~2.4]{Hoffmann}. For the converse, assume that $(f_b)_*\mu=\lambda$ for every $b\in P$. Fix $\theta(z,y)\in L(\U)$ and choose a small model $N\supseteq M$ containing its parameters. For $q\in S_y(N)$, choose $b\models q$ and put $F(q)\coloneqq\mu(\theta(f(x,b),b))$ and $H(q)\coloneqq\lambda(\theta(z,b))$. Both functions are well-defined and continuous on $S_y(N)$ by definability of $\mu$ and $\lambda$. If $q\in\supp(\nu\res N)$, then $q\res M\in\supp(\nu\res M)$ and $b\in P$, so $F(q)=((f_b)_*\mu)(\theta(z,b))=H(q)$. Since $\nu\res N$ gives measure one to its support, \eqref{eq:integral} yields
\[
 \bigl(g_*(\mu\otimes\nu)\bigr)(\theta)
 = (\mu\otimes\nu)(\theta(f(x,y),y))
 = \int_{S_y(N)}F\,d(\nu\res N)=\int_{S_y(N)}H\,d(\nu\res N)
 = (\lambda\otimes\nu)(\theta).\qedhere
\]
\end{proof}

\Needspace{8\baselineskip}
The following fact is a finite consequence of~\cite[Theorem~3.13]{CGK}, which was proved using Keisler's randomization (see also~\cite[Lemma~2.5]{Hoffmann} for a short elementary proof).

\begin{fact}\label{fact:average}
Suppose $\mu\in\measureSpace{x}{\U}$ is fim. Fix $\nu\in\measureSpace{y}{\U}$ and $\tau\in\measureSpace{xy}{\U}$. Assume that for arbitrarily large $r$ there is $\lambda_r\in\measureSpace{x_1\ldots x_r y}{\U}$ such that
\[
 \lambda_r|_{x_1,\ldots,x_r}=\mu^{(r)},\qquad
 \lambda_r|_y=\nu,\qquad
 \lambda_r|_{x_i,y}=\tau\quad(1\le i\le r).
\]
Then $\tau=\mu\otimes\nu$.
\end{fact}

A measure $\mu\in\measureSpace{\bar x}{A}$, where $\bar x\coloneqq(x_i)_{i\geq1}$, is \emph{indiscernible over $A$} if $\mu(\psi(x_{i_1},\ldots,x_{i_m}))=\mu(\psi(x_1,\ldots,x_m))$ for every $m\geq1$, every $1\leq i_1<\cdots<i_m$, and every $\psi(x_1,\ldots,x_m)\in L(A)$.
The following fact is the Keisler measure formulation of Ben Yaacov's theorem that randomization preserves NIP~\cite{BY}.

\begin{fact}[{\cite[Lemma~2.10]{HPS}}]\label{fact:indiscernible-limit}
Let $M\models T$ be any model and let $\lambda\in\measureSpace{\bar x y}{M}$. If $\lambda|_{\bar x}$ is indiscernible over $M$, then $\lambda(\varphi(x_i;y))$ converges as $i\to\infty$ for every $\varphi(x;y)\in L(M)$.
\end{fact}

\section{Pushforwards of idempotents under overlapping products}\label{sec:overlapping}

Fix a small model $M\preccurlyeq\U$, a group $G\subseteq\U^x$ type-definable over $M$, and an idempotent measure $\mu\in\measureSpace{G}{\U}$ generically stable over $M$.

\begin{definition}\label{def:word-map}
Fix $M$-definable maps $\widetilde m\colon(\U^x)^2\to\U^x$ and $\widetilde\iota\colon\U^x\to\U^x$ extending multiplication and inversion on $G$~\cite[\S2.2]{CGK}. A group word $w(x_1,\ldots,x_n)$ defines a map $w\colon G^n\to G$. Its $M$-definable extension $\widehat w\colon(\U^x)^n\to\U^x$ interprets the variables $x_i$ as coordinate projections, the identity as $1_G$, and multiplication and inversion as $\widetilde m$ and $\widetilde\iota$, evaluating products from left to right. Write $\bar u\coloneqq(x_1,\ldots,x_r)$ and $\bar v\coloneqq(x_{r+1},\ldots,x_n)$. For $\bar b\in G^{n-r}$, the map $\bar a\mapsto\widehat w(\bar a,\bar b)$ is an $M\bar b$-definable extension of $w(\bar u,\bar b)$.

Throughout the rest of the paper we write $w$ for $\widehat w$ when using a word map. In particular, $w_*\lambda$ means $(\widehat w)_*\lambda$.
\end{definition}

\begin{remark}\label{rem:agreement}
Let $\lambda\in\measureSpace{x_1\ldots x_n}{\U}$ with $\lambda|_{x_i}\in\measureSpace{G}{\U}$ for every $i$. Then $\lambda$ concentrates on $G^n$. If definable maps $f,g\colon(\U^x)^n\to\U^z$ agree on $G^n$, the definable set $\{f=g\}$ contains $G^n$, so $\lambda(f=g)=1$ and $f_*\lambda=g_*\lambda$. In particular, two word maps that agree on $G^n$ have the same pushforward of $\lambda$.
\end{remark}

\begin{definition}
For $b\in G$ and $m\geq1$, consider the word map $w_{m,b}\colon(x_0,\ldots,x_m)\mapsto x_0 b x_1\cdots b x_m$. Put $w_{0,b}\colon x_0\mapsto x_0$.
\end{definition}

The following lemma generalizes~\cite[Lemma~3.1]{Hoffmann} (which is the fim measure analogue of~\cite[Remark~2.22]{CGK} for invariant idempotent types) to type-definable groups.

\begin{lemma}\label{lem:sandwich}
For every $b\models p\in\supp(\mu\res M)$,
\begin{equation}\label{eq:fixed}
      \mu*\delta_b*\mu=\mu .
\end{equation}
\end{lemma}
\begin{proof}
Let $S\colon(u,y,v)\mapsto(y,uyv)$ and put $\tau\coloneqq S_*\mu^{(3)}$. For $n\geq1$, write $\bar x\coloneqq(x_0,\ldots,x_{n+1})$ and define $H_n\colon\bar x\mapsto(x_1,\ldots,x_n,x_0\cdots x_{n+1})$. Let $\lambda_n\coloneqq(H_n)_*\mu^{(n+2)}$ in target variables $(x_1,\ldots,x_n,z)$. Then $\lambda_n|_{x_1,\ldots,x_n}=\mu^{(n)}$, and $\lambda_n|_z=\mu$ by Remark~\ref{rem:blocks}\textup{(\ref{part:blocks-one})}.

For $1\leq i\leq n$, the map $B_i\colon\bar x\mapsto(x_0\cdots x_{i-1},x_i,x_{i+1}\cdots x_{n+1})$ sends $\mu^{(n+2)}$ to $\mu^{(3)}$ by Remark~\ref{rem:blocks}\textup{(\ref{part:blocks-several})}. For $\bar a=(a_0,\ldots,a_{n+1})\in G^{n+2}$ we have $(\pi_{x_i,z}\circ H_n)(\bar a)=(a_i,a_0\cdots a_{n+1})=(S\circ B_i)(\bar a)$. Hence Remark~\ref{rem:agreement} gives $\lambda_n|_{x_i,z}=S_*((B_i)_*\mu^{(n+2)})=\tau$. Fact~\ref{fact:average} now gives $\tau=\mu^{(2)}$, so $S_*\mu^{(3)}=\mu^{(2)}$. Let $f(u,v;y)\coloneqq uyv$. By symmetry,
\begin{equation}\label{eq:one-overlap}
 \bigl((u,v,y)\mapsto(f(u,v;y),y)\bigr)_*
       (\mu^{(2)}_{u,v}\otimes\mu_y)=\mu_z\otimes\mu_y .
\end{equation}

Fix $b\models p\in\supp(\mu\res M)$. Then $b\in G$ since $G$ is type-definable over $M$ and $\mu$ concentrates on $G$. Applying the forward direction of Lemma~\ref{lem:fibres} to~\eqref{eq:one-overlap} gives $\bigl((u,v)\mapsto ubv\bigr)_*\mu^{(2)}_{u,v}=\mu_z$.

Let $R_b(u)\coloneqq ub$ and $m(s,v)\coloneqq sv$. On $G^2$ we have $w_{1,b}=m\circ(R_b\times\mathrm{id})$. Since $(R_b)_*\mu=\mu*\delta_b$, Fact~\ref{fact:products}\textup{(\ref{part:pushforwards})} gives $(R_b\times\mathrm{id})_*\mu^{(2)}=(\mu*\delta_b)\otimes\mu$. Hence $\mu=(w_{1,b})_*\mu^{(2)}=m_*\bigl((\mu*\delta_b)\otimes\mu\bigr)=(\mu*\delta_b)*\mu$ by Remark~\ref{rem:agreement} and the definition of convolution.
\end{proof}

\begin{corollary}\label{cor:repeat}
If $b\models p\in\supp(\mu\res M)$, then $(w_{m,b})_*\mu^{(m+1)}=\mu$ for every $m\geq0$.
\end{corollary}
\begin{proof}
The map $w_{0,b}$ is the identity. For the induction step, write
\[
 w_{m+1,b}(x_0,\ldots,x_{m+1})
 =w_{1,b}\bigl(w_{m,b}(x_0,\ldots,x_m),x_{m+1}\bigr).
\]
By the induction hypothesis and Remark~\ref{rem:blocks}\textup{(\ref{part:blocks-general})}, the inner pair of maps sends $\mu^{(m+2)}$ to $\mu^{(2)}$. Applying $w_{1,b}$ gives $\mu$ by Lemma~\ref{lem:sandwich}.
\end{proof}

The following is the key technical result.

\begin{lemma}\label{lem:products}
For $k\ge1$, let $\bar x\coloneqq(x_0,\ldots,x_k)$ and $\bar y\coloneqq(y_1,\ldots,y_k)$ be tuples of variables. Consider the $M$-definable maps $F_k$ and $f_k$ on $\U^{\bar x\bar y}$ with values in $(\U^x)^k$ and $(\U^x)^{2k}$, respectively, given by
\[
     F_k(\bar x,\bar y)\coloneqq (x_{j-1}y_jx_j:1\leq j\leq k),
     \qquad f_k(\bar x,\bar y)\coloneqq (F_k(\bar x,\bar y),\bar y).
\]
Then
\begin{equation}\label{eq:products}
           (f_k)_*(\mu^{(2k+1)})=\mu^{(k)}\otimes\mu^{(k)} .
\end{equation}
\end{lemma}
\begin{proof}
For a global Keisler measure $\rho$ in variables $w$ and a tuple $g=(g_1,\ldots,g_m)$ of definable maps on $\U^w$, write $\mathcal L_\rho(g_1,\ldots,g_m)\coloneqq g_*\rho$. Thus $\mathcal L_\rho(g)(\varphi)=\rho(\varphi(g(w)))$ for every $L(\U)$-formula $\varphi$ in the target variables. When a subtuple $v$ of the source variables $w$ is listed among maps, it denotes the corresponding coordinate projection map $\pi_v\colon\U^w\to\U^v$.

We prove~\eqref{eq:products} for every $k\geq1$ by induction on $k$. For $k=1$ this is~\eqref{eq:one-overlap}. Fix $k\geq2$, assume~\eqref{eq:products} for $k-1$, and put $\nu\coloneqq\mu^{(k-1)}$. We must show that $\tau\coloneqq(f_k)_*\mu^{(2k+1)}$ equals $\mu^{(k)}\otimes\mu^{(k)}$.

Let $\bar x'\coloneqq(x_0,\ldots,x_{k-1})$ and $\bar y'\coloneqq(y_1,\ldots,y_{k-1})$. By the definition of $F_k$,
\begin{equation}\label{eq:last-product}
 F_k(\bar x,\bar y)=\bigl(F_{k-1}(\bar x',\bar y'),\,x_{k-1}y_kx_k\bigr).
\end{equation}
We derive a consequence of the induction hypothesis for use in Claim~\ref{cl:rows}. Let $\bar s$ be an additional tuple of variables. For every global measure $\gamma$ in $(\bar y',\bar s)$ with $\gamma|_{\bar y'}=\nu$, we have
\begin{equation}\label{eq:retain-outputs}
 \mathcal L_{\mu^{(k)}_{\bar x'}\otimes\gamma}\bigl(F_{k-1}(\bar x',\bar y'),\bar s\bigr)=\nu\otimes(\gamma|_{\bar s}).
\end{equation}
Indeed, $\mu^{(2k-1)}=\mu^{(k)}_{\bar x'}\otimes\nu_{\bar y'}$ by associativity and symmetry of the Morley powers. Thus the induction hypothesis $(f_{k-1})_*\mu^{(2k-1)}=\nu\otimes\nu$ says that $\mathcal L_{\mu^{(k)}_{\bar x'}\otimes\nu_{\bar y'}}(F_{k-1}(\bar x',\bar y'),\bar y')=\nu\otimes\nu$. By the forward direction of Lemma~\ref{lem:fibres} with parameter variables $\bar y'$, the map $\bar a\mapsto F_{k-1}(\bar a,\bar b)$ sends $\mu^{(k)}$ to $\nu$ whenever $\tp(\bar b/M)\in\supp(\nu\res M)$. Now let $(\bar b,\bar c)$ realize a type in $\supp(\gamma\res M)$. Since $\gamma|_{\bar y'}=\nu$, we have $\tp(\bar b/M)\in\supp(\nu\res M)$. The fibre of $(\bar x',\bar y',\bar s)\mapsto F_{k-1}(\bar x',\bar y')$ at $(\bar y',\bar s)=(\bar b,\bar c)$ is precisely $\bar a\mapsto F_{k-1}(\bar a,\bar b)$. The reverse direction of Lemma~\ref{lem:fibres} therefore gives
\[
 \mathcal L_{\mu^{(k)}_{\bar x'}\otimes\gamma}
   \bigl(F_{k-1}(\bar x',\bar y'),\bar y',\bar s\bigr)=\nu\otimes\gamma.
\]
Push both measures in this equality forward by $(\bar v,\bar y',\bar s)\mapsto(\bar v,\bar s)$, where $\bar v$ is the tuple of output variables of $F_{k-1}$. By Fact~\ref{fact:products}\textup{(\ref{part:pushforwards})}, this sends $\nu\otimes\gamma$ to $\nu\otimes(\gamma|_{\bar s})$, proving~\eqref{eq:retain-outputs}.

Fix $r\geq1$. Introduce tuples of variables $\bar z\coloneqq(z_0,\ldots,z_{r+1})$ and $\bar u_i\coloneqq(u_{i,0},\ldots,u_{i,k-2})$ for $1\leq i\leq r$. Let $\omega\coloneqq\mu^{(rk+2)}_{\bar z,\bar u_1,\ldots,\bar u_r}\otimes\mu^{(k)}_{\bar y}$. Write $\mathcal L\coloneqq\mathcal L_\omega$. In order to apply Fact~\ref{fact:average}, we construct maps $t_1,\ldots,t_r,V$ with $\mathcal L(t_1,\ldots,t_r)=\nu^{(r)}$, $\mathcal L(V)=\mu^{(k+1)}$, and $\mathcal L(t_i,V)=\tau$ for every $i$. To obtain the same $V$ for each $i$, we use $r$ factorizations of the word $P_{r+1}$ below. Define
\[
\begin{aligned}
 P_j(\bar z,y_k)&\coloneqq z_0y_kz_1\cdots y_kz_j &&(0\leq j\leq r+1),\\
 S_i(\bar z,y_k)&\coloneqq z_{i+1}y_k\cdots y_kz_{r+1} &&(1\leq i\leq r),\\
 t_i(\bar z,\bar u_i,\bar y)&\coloneqq
 \bigl((u_{i,j-1}y_j u_{i,j}:1\leq j\leq k-2),\ u_{i,k-2}y_{k-1}P_i\bigr) &&(1\leq i\leq r).
\end{aligned}
\]
Here $P_0=z_0$, $S_r=z_{r+1}$, and for $k=2$ the map $t_i$ has just one coordinate. Put $V\coloneqq(P_{r+1},\bar y)$.

\begin{claim}\label{cl:shared}
$\mathcal L(P_j,\bar y)=\mu^{(k+1)}$ for $0\leq j\leq r+1$, and $\mathcal L(t_i,V)=\tau$ for $1\leq i\leq r$.
\end{claim}
\begin{proof}[Proof of Claim~\ref{cl:shared}]
Fix a tuple of elements $\bar b=(b_1,\ldots,b_k)$ realizing a type in $\supp(\mu^{(k)}\res M)$. Its last coordinate $b_k$ realizes a type in $\supp(\mu\res M)$. The maps $P_j(\bar z,b_k)$ and $S_i(\bar z,b_k)$ each send $\mu^{(rk+2)}$ to $\mu$ by Corollary~\ref{cor:repeat} applied to $(z_0,\ldots,z_j)$ and $(z_{i+1},\ldots,z_{r+1})$, respectively. The three blocks of variables $\bar u_i$, $(z_0,\ldots,z_i)$ and $(z_{i+1},\ldots,z_{r+1})$ are disjoint. Hence the map $(\bar u_i,P_i(\bar z,b_k),S_i(\bar z,b_k))$ sends $\mu^{(rk+2)}$ to $\mu^{(k+1)}$ by Remark~\ref{rem:blocks}\textup{(\ref{part:blocks-general})}. The reverse direction of Lemma~\ref{lem:fibres}, with parameter variables $\bar y$, gives
\[
 \mathcal L(P_j,\bar y)=\mu^{(k+1)},\qquad
 \mathcal L(\bar u_i,P_i,S_i,\bar y)=\mu^{(2k+1)}.
\]
Substituting $(\bar u_i,P_i,S_i)$ for $\bar x$ in~\eqref{eq:last-product} gives
\[
 F_k\bigl((\bar u_i,P_i,S_i),\bar y\bigr)
 =\bigl(F_{k-1}((\bar u_i,P_i),\bar y'),P_i y_k S_i\bigr).
\]
The $j$-th output of $F_{k-1}((\bar u_i,P_i),\bar y')$ is $u_{i,j-1}y_j u_{i,j}$ for $1\le j\le k-2$, and its last output, in position $k-1$, is $u_{i,k-2}y_{k-1}P_i$. These are the coordinates defining $t_i$. On $G^{rk+k+2}$ we have $P_i y_k S_i=P_{r+1}$, so $f_k\circ(\bar u_i,P_i,S_i,\bar y)=(t_i,V)$ there. Since $\omega$ concentrates on $G^{rk+k+2}$, Remark~\ref{rem:agreement} gives equal pushforwards of $\omega$ under these two maps. Composing pushforwards and using the value of $\mathcal L(\bar u_i,P_i,S_i,\bar y)$ computed above gives $\mathcal L(t_i,V)=(f_k)_*\mathcal L(\bar u_i,P_i,S_i,\bar y)=(f_k)_*\mu^{(2k+1)}=\tau$.
\end{proof}

\begin{claim}\label{cl:rows}
$\mathcal L(t_1,\ldots,t_r)=\nu^{(r)}$.
\end{claim}
\begin{proof}[Proof of Claim~\ref{cl:rows}]
Fix $1\le i\le r$. Let $U_i\colon(\U^x)^{rk+k+2}\to(\U^x)^k$ be the coordinate projection $(\bar z,\bar u_1,\ldots,\bar u_r,\bar y)\mapsto(\bar u_i,z_i)$, and let $D_i$ be the map on the same domain with coordinates $(y_1,\ldots,y_{k-2},y_{k-1}P_{i-1}y_k)$. By Claim~\ref{cl:shared} with $j=i-1$, $\mathcal L(P_{i-1},\bar y)=\mu^{(k+1)}$. Keeping $y_1,\ldots,y_{k-2}$ and multiplying $y_{k-1},P_{i-1},y_k$ in that order gives $\mathcal L(D_i)=\mu^{(k-1)}=\nu$ by Remark~\ref{rem:blocks}\textup{(\ref{part:blocks-several})}.

Let $T_{i-1}\coloneqq(t_1,\ldots,t_{i-1})$, with $T_0$ empty, and put $\gamma\coloneqq\mathcal L(D_i,T_{i-1})$ in target variables $(\bar y',\bar s)$. The maps $D_i$ and $T_{i-1}$ depend only on $\bar y$, $z_0,\ldots,z_{i-1}$, and $\bar u_1,\ldots,\bar u_{i-1}$. These variables are disjoint from $(\bar u_i,z_i)$, so $\mathcal L(U_i,D_i,T_{i-1})=\mu^{(k)}\otimes\gamma$ by Remark~\ref{rem:blocks}\textup{(\ref{part:blocks-general})}. Since $\gamma|_{\bar y'}=\mathcal L(D_i)=\nu$, \eqref{eq:retain-outputs} says that $(\bar x',\bar y',\bar s)\mapsto(F_{k-1}(\bar x',\bar y'),\bar s)$ sends $\mu^{(k)}\otimes\gamma$ to $\nu\otimes(\gamma|_{\bar s})$. Its composite with $(U_i,D_i,T_{i-1})$ is $(F_{k-1}\circ(U_i,D_i),T_{i-1})$ and therefore sends $\omega$ to $\nu\otimes(\gamma|_{\bar s})$.

For $1\le j\le k-2$, the $j$-th output of $F_{k-1}\circ(U_i,D_i)$ is $u_{i,j-1}y_j u_{i,j}$, the $j$-th output of $t_i$. On $G^{rk+k+2}$ its $(k-1)$-st output is $u_{i,k-2}(y_{k-1}P_{i-1}y_k)z_i=u_{i,k-2}y_{k-1}P_i$, the last output of $t_i$, since $P_i=P_{i-1}y_kz_i$. By Remark~\ref{rem:agreement}, the maps $(F_{k-1}\circ(U_i,D_i),T_{i-1})$ and $(t_i,T_{i-1})$ thus have the same pushforward of $\omega$. The $\bar s$-coordinates of $(D_i,T_{i-1})$ are $T_{i-1}$, so $\mathcal L(t_i,T_{i-1})=\nu\otimes(\gamma|_{\bar s})=\nu\otimes\mathcal L(T_{i-1})$. For $i=1$ this gives $\mathcal L(t_1)=\nu$. Induction on $i$ now gives $\mathcal L(t_1,\ldots,t_i)=\nu^{(i)}$ by the definition of the Morley powers.
\end{proof}

Take $\lambda_r\coloneqq\mathcal L(t_1,\ldots,t_r,V)$. The restriction of $\lambda_r$ to its last $k+1$ coordinates is $\mathcal L(V)=\mathcal L(P_{r+1},\bar y)=\mu^{(k+1)}$ by the definition of $V$ and the first assertion of Claim~\ref{cl:shared} with $j=r+1$. The restriction of $\lambda_r$ to the coordinates of $(t_1,\ldots,t_r)$ is $\nu^{(r)}$ by Claim~\ref{cl:rows}, and its restriction to the coordinates of $(t_i,V)$ is $\tau$ for each $1\le i\le r$ by the second assertion of Claim~\ref{cl:shared}. Since $\nu$ is fim by Fact~\ref{fact:products}\textup{(\ref{part:fimclosure})}, Fact~\ref{fact:average} gives $\tau=\nu\otimes\mu^{(k+1)}=\mu^{(k)}\otimes\mu^{(k)}$.
\end{proof}

\section{Invariance under inversion}

\begin{proposition}\label{prop:sym}
Every generically stable idempotent $\mu\in\measureSpace{G}{\U}$ is invariant under inversion, i.e., $\mu=\mu^{-1}$, where $\mu^{-1}\coloneqq(x\mapsto x^{-1})_*\mu$.
\end{proposition}

\begin{proof}
Fix $\alpha(x)\in L(\U)$. Enlarge $M$ so that $\alpha\in L(M)$ and $\mu$ is generically stable over $M$. The group $G\times G$ acts on $G$ on the right by $t\cdot(u,v)\coloneqq u^{-1}tv$. Let $\omega\coloneqq\mu^{(\omega)}$ in the variables $(x_0,y_1,x_1,y_2,x_2,\ldots)$. For $e=(a_0,b_1,a_1,b_2,a_2,\ldots)\in(\U^x)^\omega$, set $c_i\coloneqq a_{i-1}b_i a_i$ and define $h_i\colon(\U^x)^\omega\to(\U^x)^2$ for $i\geq1$ by
\[
 h_i(e)\coloneqq
 \begin{cases}
   (c_i,b_i) & i\text{ odd},\\
   (b_i,c_i) & i\text{ even}.
 \end{cases}
\]
Define $Z_i\colon(\U^x)^\omega\to(\U^x)^2$ by $Z_i(e)\coloneqq h_1(e)\cdots h_i(e)$, multiplying pairs componentwise. We prove by induction that whenever the first $2i+1$ entries of $e$ lie in $G$, $a_0\cdot Z_i(e)$ is $a_i$ for even $i$ and $a_i^{-1}$ for odd $i$. With $Z_0(e)\coloneqq(1_G,1_G)$, we have $a_0\cdot Z_0(e)=a_0$. For $i\geq1$, the right action gives $a_0\cdot Z_i(e)=(a_0\cdot Z_{i-1}(e))\cdot h_i(e)$. Substituting the induction hypothesis and using $c_i=a_{i-1}b_i a_i$ gives
\begin{equation}\label{eq:alternating-action}
 a_0\cdot Z_i(e)=
 \begin{cases}
   a_{i-1}\cdot(c_i,b_i)=c_i^{-1}a_{i-1}b_i=a_i^{-1} & i\text{ odd},\\
   a_{i-1}^{-1}\cdot(b_i,c_i)=b_i^{-1}a_{i-1}^{-1}c_i=a_i & i\text{ even}.
 \end{cases}
\end{equation}

Write $\bar z\coloneqq(z_i)_{i\geq1}$, where each $z_i$ is a pair of $x$-variables. Let $\lambda\in\measureSpace{\bar z t}{M}$ be the restriction to $M$ of $(e\mapsto((Z_i(e))_{i\geq1},a_0))_*\omega$. We show that $\lambda|_{\bar z}$ is indiscernible over $M$. For $n\geq1$, let $\pi_n(e)\coloneqq(a_0,\ldots,a_n,b_1,\ldots,b_n)$. This reorders the first $2n+1$ coordinates of $e$, so $(\pi_n)_*\omega=\mu^{(2n+1)}$ by symmetry. Evaluate the map $f_n$ of Lemma~\ref{lem:products} at $\pi_n(e)$, taking $x_j=a_j$ for $0\leq j\leq n$ and $y_j=b_j$ for $1\leq j\leq n$. Its first $n$ outputs are $a_{j-1}b_j a_j=c_j$ in order, and its last $n$ outputs are $b_1,\ldots,b_n$. Thus~\eqref{eq:products} gives
\[
 (e\mapsto(c_1,\ldots,c_n,b_1,\ldots,b_n))_*\omega
 =(f_n)_*((\pi_n)_*\omega)=\mu^{(n)}\otimes\mu^{(n)}=\mu^{(2n)}.
\]
Reordering $(c_1,\ldots,c_n,b_1,\ldots,b_n)$ into the pairs $h_1(e),\ldots,h_n(e)$ preserves $\mu^{(2n)}$ by Fact~\ref{fact:products}\textup{(\ref{part:symmetry})}, so $(h_1,\ldots,h_n)_*\omega=\mu^{(2n)}$.

Given $1\leq i_1<\cdots<i_m$, put $i_0\coloneqq0$ and $N\coloneqq i_m$. Define $D_j\colon(\U^x)^\omega\to(\U^x)^2$ for $1\leq j\leq m$ by $D_j(e)\coloneqq h_{i_{j-1}+1}(e)\cdots h_{i_j}(e)$, multiplying the pairs componentwise in increasing order of their indices. Thus the $2m$ coordinates of $(D_1,\ldots,D_m)$ are products on disjoint blocks of the $2N$ coordinates of $(h_1,\ldots,h_N)$. By Remark~\ref{rem:blocks}\textup{(\ref{part:blocks-several})},
\begin{equation}\label{eq:block-increments}
 (D_1,\ldots,D_m)_*\omega=\mu^{(2m)}=(h_1,\ldots,h_m)_*\omega.
\end{equation}
Push~\eqref{eq:block-increments} forward under the map $C_m\colon(g_1,\ldots,g_m)\mapsto(g_1,g_1g_2,\ldots,g_1\cdots g_m)$, multiplying pairs componentwise. Both compositions depend only on the first $2N+1$ entries of $e$. On $G^{2N+1}$ they satisfy
\[
 C_m\circ(D_1,\ldots,D_m)=(Z_{i_1},\ldots,Z_{i_m}),\qquad
 C_m\circ(h_1,\ldots,h_m)=(Z_1,\ldots,Z_m).
\]
Indeed, $D_1(e)\cdots D_j(e)=Z_{i_j}(e)$ because the first $j$ intervals partition $1,\ldots,i_j$ in order, and $h_1(e)\cdots h_j(e)=Z_j(e)$ by definition. The restriction of $\omega$ to these coordinates concentrates on $G^{2N+1}$, so Remark~\ref{rem:agreement} gives
$(Z_{i_1},\ldots,Z_{i_m})_*\omega=(Z_1,\ldots,Z_m)_*\omega$. This proves indiscernibility.

By Fact~\ref{fact:indiscernible-limit}, $\lambda(\alpha(t\cdot z_i))$ converges as $i\to\infty$. By the definition of $\lambda$, its $i$-th term is $\bigl((e\mapsto a_0\cdot Z_i(e))_*\omega\bigr)(\alpha)$. On $G^{2i+1}$, \eqref{eq:alternating-action} identifies this map with $e\mapsto a_i^{-1}$ for odd $i$ and with $e\mapsto a_i$ for even $i$. Hence by Remark~\ref{rem:agreement} and $\omega|_{x_i}=\mu$, the odd terms are $\mu(\alpha(x^{-1}))$ and the even terms are $\mu(\alpha(x))$. Their common limit proves the assertion.
\end{proof}

\section{Concentration on the stabilizer}

\begin{proof}[Proof of Theorem~\ref{thm:main}]
Choose a small model $M$ over which $G$ is type-definable and $\mu\in\measureSpace{G}{\U}$ is generically stable. First suppose $\mu$ is idempotent. Fix $b\models p\in\supp(\mu\res M)$. Then $b\in G$ and $\mu*\delta_b*\mu=\mu$ by Lemma~\ref{lem:sandwich}. Set $\nu\coloneqq\delta_b*\mu=(x\mapsto bx)_*\mu$. Since $x\mapsto bx$ is $Mb$-definable, $\nu$ is generically stable over any small model containing $Mb$ (by Facts~\ref{fact:fim} and~\ref{fact:products}\textup{(\ref{part:fimclosure})}). Associativity gives
\[
       \nu*\nu=\delta_b*(\mu*\delta_b*\mu)=\delta_b*\mu=\nu .
\]
By Proposition~\ref{prop:sym}, $\mu=\mu^{-1}$ and $\nu=\nu^{-1}$. Since $(bx)^{-1}=x^{-1}b^{-1}$ on $G$, Remark~\ref{rem:agreement} and composition of pushforwards give
\[
\begin{aligned}
 \delta_b*\mu
 &=\nu=\nu^{-1}
   =(y\mapsto y^{-1})_*\bigl((x\mapsto bx)_*\mu\bigr)
   =(x\mapsto(bx)^{-1})_*\mu\\
 &=(x\mapsto x^{-1}b^{-1})_*\mu
   =(y\mapsto yb^{-1})_*\mu^{-1}
   =\mu^{-1}*\delta_{b^{-1}}=\mu*\delta_{b^{-1}}.
\end{aligned}
\]
Convolving on the left by $\mu$ and using \eqref{eq:fixed} and idempotence yields $\mu=\mu*\delta_{b^{-1}}$. Thus $b^{-1}$ and hence $b$ belong to $H=\Stab_{\mathrm r}(\mu)$. By \cite[Fact~3.35]{CGK}, $H$ is type-definable over $M$. Every realization of every type in $\supp(\mu\res M)$ belongs to $H$ by the preceding argument. Therefore the clopen set corresponding to each formula in an $M$-definition of $H$ contains $\supp(\mu\res M)$ and has measure one. This proves that $\mu$ concentrates on $H$. By \cite[Remark~3.39]{CGK}, $H=\Stab_{\ell}(\mu)$ and $H$ is a fim group. Among global measures concentrating on $H$, $\mu$ is the unique left-$H$-invariant measure and the unique right-$H$-invariant measure. Since $H$ is fim, it is fsg by \cite[Remark~3.30(1)]{CGK}.

Conversely, if $\mu$ concentrates on its right stabilizer $H$, then it is a right-invariant measure on $H$ and $\mu$ is idempotent (see \cite[Proposition~3.19]{CG}).
\end{proof}

\section*{Acknowledgements}
We thank Kyle Gannon, Krzysztof Krupi\'nski, Sergei Starchenko, and Atticus Stonestrom for helpful conversations and collaborations on this topic over the years. Chernikov was partially supported by the NSF Research Grants DMS-2246598, DMS-2554164 and by the Alexander von Humboldt Foundation.

\section*{AI methodology}
The result was established over several months of research with various AI models: ChatGPT 5.6 and 6; Claude Opus 5.5 and Fable 5.2. In particular, the solvable case, with some further extensions, was established in early August 2026 using GPT5.6 (with some overlap of the proof with the one in~\cite{Hoffmann}). The other ingredients leading to the case of types in strongly dependent theories were developed by early November 2026 using GPT6, with the final key technical result (full Lemma~\ref{lem:products}) obtained in mid November. The original proof was significantly reworked by the author, in parts with further AI assistance (using Codex, GPT6.2 and Claude Opus 5.5).

\printbibliography
\end{document}